\documentclass[a4paper,10pt]{article}
\usepackage{palatino}
\usepackage{graphicx}
\usepackage{graphics,color}
\usepackage{amsmath}
\usepackage{amsopn}
\usepackage{amstext}
\usepackage{rotating}
\usepackage{amsthm}
\usepackage{amsbsy}
\usepackage{amsfonts}
\usepackage{amssymb}
\usepackage{layout}
\usepackage{latexsym,lscape, graphicx, float}
\begin{document}

\title{The Bismut-Elworthy-Li formula for jump-diffusions and applications
to Monte Carlo methods in finance} \author{Thomas R. Cass and Peter K.
Friz\\Statistical Laboratory, University of Cambridge,\\Wilberforce Road,
Cambridge, CB3 0WB, UK.} \date{March 22, 2007} \maketitle

\begin{abstract} We extend the Bismut-Elworthy-Li formula to non-degenerate
jump diffusions and "payoff" functions depending on the process at multiple
future times. In the spirit of Fourni\'e et al $[14]$ and Davis and Johansson
$[10]$ this can improve Monte Carlo numerics for stochastic volatility models
with jumps.
To this end one needs so-called Malliavin weights and we give explicit formulae
valid
in presence of jumps: (a) In a non-degenerate situation, the extended BEL
formula represents
possible Malliavin weights as Ito integrals with explicit integrands;
(b) in a hypoelliptic setting
we review work of Arnaudon and Thalmaier $[1]$ and also find explicit weights, now
involving the Malliavin covariance matrix, but still straight-forward to implement. (This is
in contrast to recent work by
Forster, L\"utkebohmert and Teichmann where weights are constructed as
anticipating Skorohod integrals.)
We give some financial examples covered by (b) but note that most practical
cases of poor Monte Carlo performance,
Digital Cliquet contracts for instance, can be dealt with by the extended BEL
formula and hence without any reliance
on Malliavin calculus at all. We then discuss some of the approximations, often
ignored in the literature,
needed to justify the use of the Malliavin weights in the context of standard
jump diffusion models. Finally, as all this
is meant to improve numerics, we give some numerical results with focus on
Cliquets under the Heston model with jumps.
\end{abstract}

\newtheorem{axiom}{Lemma}

\newtheorem{theorem}{Theorem}%

\newtheorem{subaxiom}{Corollary}

\newtheorem{subtheorem}{Definition}

\newtheorem{subsubtheorem}{Proposition}

\newtheorem{subsubaxiom}{Remark}%

\newtheorem{subsubsubtheorem}{Assumption}

\newtheorem{subsubsubaxiom}{Example}

\setcounter{section}{0}

\newpage

\section{Introduction}

Modern arbitrage theory reduces the pricing of (non-American) options to
the computation of an expectation under a risk neutral measure. It is
common practice to assume that the risk neutral measure is induced by a
parametric family of jump diffusions which can then be calibrated to liquid
option prices. We can therefore assume that all expectations are with
respect to a fixed pricing measure. A typical option on some underlying
$(S_{t})$ then has (undiscounted) price \[ \mathbb{E}\left[ f\left(
S_{T_{1}},S_{T_{2}},\ldots,S_{T_{n}}\right) \right]
\equiv\mathbb{E}[f(\underline{S})]. \] For hedging and risk-management
purposes it is crucial to understand the dependence on $S_{0}$ and other
model parameters. Computing \[ \Delta=\frac{\partial}{\partial
S_{0}}\mathbb{E}\left[ f\left( \underline {S}\right) \right]
=\mathbb{E}\left[ \nabla f\left( \underline{S}\right)
\frac{\partial\underline{S}}{\partial S_{0}}\right] \] via finite
differences can present computational challenges in Monte Carlo; just think
of an at-the-money digital option near expiration. Broadie and Glasserman
$[7]$ showed that this problem is overcome by \begin{align}
 \frac{\partial}{\partial
S_{0}}\mathbb{E}\left[ f\left( \underline {S}\right) \right]
=\mathbb{E}\left[ f\left( \underline{S}\right) \pi\right]
\end{align} where $\pi$ is
the logarithmic derivative of the joint density of the random vector
$\underline{S}$. On the other hand, the random weight $\pi$ adds noise
itself and it is important to localise: for instance by using $(1)$ for an
irregular, but compactly support and bounded, $\tilde{f}$ and the usual
finite difference technique for $f-\tilde{f}$, assumed to be nice ($C^{1}$
will usually suffice).

In two seminal papers, Fourni\'e et al $[14]$ and $[15]$ use Malliavin
calculus to compute $\pi$ when no explicit transition density is known.
They work with non-degenerate (or: elliptic) continuous diffusions but also
cover some hypoelliptic situations. As is well known, elliptic results can
be obtained by the Bismut-Elworthy-Li formula (Elworthy and Li $[11]$,
Bismut $[6]$ ) and there are, in fact, other ways to obtain such results
without Malliavin calculus: we mention in particular the idea of Thalmaier
$[27]$ of differentiation at the level of local martingales which was employed by Gobet and Munos
$[18]$ in the present context. The point was that in many cases of practical interest,
at least in absence of jumps, one does not need Malliavin calculus. (Specialists will
note that Malliavin techniques are more flexible in the sense that different
perturbations of Brownian motion yield different weights and there is an
apriori interest to pick weights with small variance. In reality, it is hard to
justify much effort in this direction as the potential gains are
negligible to the improvements obtained by localisation.)

Over the last decade it has become clear that pure diffusion models are
unable to fit the short-dated smile and jumps have been included to models
to rectify this situation; Cont and Tankov $[9]$ and Gatheral $[16]$
provide two excellent accounts. The question has arisen as to how the above
ideas can be adapted to models based on jump diffusion processes and we
shall propose a quite simple solution to this along the ideas of
Elworthy-Li bypassing both classical Malliavin techniques and its
extensions to L\'{e}vy processes that have been used in this financial Monte
Carlo context.
We note that a similar extension of the BEL formula, slightly less general than
ours, was used recently
by Priola and Zabczyk $[22]$ to establish Lioville theorems for non-local
operators.

\bigskip

Let us briefly mention that in some cases a random weight $\pi$ can be
constructed by conditioning arguments. Consider for instance the trivial
example $X_{t}=z+B_{t}+N_{t}$, where $B$ is a standard Brownian motion and
$N$ a Poisson process. Conditional on $N_{t}$, any function of $X_{t}$ is a
(different) function of $z+B_{t}$, a pure diffusion with no jumps, and
since the associated random weight $\pi$ is universal\ (i.e. do not depend
on the particular payoff function) this also solves the problem for the
jump diffusion $X$. This kind of reasoning leads immediately to the class
of "separable" jump diffusions, considered in Davis and Johansson $[10]$ via
Malliavin calculus for simple L\'{e}vy processes. We shall omit a detailed
discussion since a refined, iterated conditioning argument can be used
assuming {\it only} finite activity of the jumps (and {\it without} assuming
separability in the sense of $[10]$). To this
end, we quickly recall the BEL for continuous diffusions (see Section 3 for
notation and assumptions) \[ \frac{\partial}{\partial
z_{j}}\mathbb{E}[f(x_{T}^{z})]=\mathbb{E}\left[
f(x_{T}^{z})\int_{0}^{T}a\left( t\right) \left( R(t,x_{t}^{z}%
)\frac{\partial x_{t}^{z}}{\partial z_{j}}\right) ^{T}dW_{t}\right] \]
where $\int_{0}^{T}a\left( t\right) dt=1$ and $x_{0}=z$. Let $0<S<T$ be
deterministic. The standard choice $a\equiv1/T$ gives a weight, say $\pi
_{0,T}$. Another weight (of higher variance) comes from
\thinspace$a\equiv1/S$ on $\left[ 0,S\right] ,$ $0$ otherwise, and we call
it $\pi_{0,S}$. One can also condition on $x_{S}$ and apply the BEL formula
over the time interval $\left[ S,T\right] $, this yields another weight
$\pi_{S,T}$ for the derivative of $\mathbb{E}[f(x_{T})|x_{S}]$
w.r.t. $x_{S}$. We leave it to the reader to check that,
combined with the chain-rule, $\left( \partial/\partial z_{j}\right)
=\left( \partial x_{S}/\partial z\right) \partial/\partial x_{S}$, the
weight $\pi_{0,T}$ can be assembled from $\pi_{0,S}$ and $\pi_{S,T}$. In
other words, instead of applying BEL on $[0,T)$ one can apply it on $[0,S)$
and $[S,T)$. While we did not assume jumps in this discussion, it is clear
that a cadlag discontinuity of $x$ at time $S$ does not pose a problem.
This extends to any number of intervals and if we are dealing with a finite
activity jump diffusion conditioning will reduce the problem to the one
just discussed. \textit{The flaw with this sort of reasoning is that it
makes fundamental use of a property which is completely irrelevant for the
result to hold true: finite activity of jumps}. In the general case, i.e.
beyond
finite activity, not only does the preceding argument break down, but jumps
arise from a genuine stochastic integral w.r.t. a compensated Poisson
random measure and any conditioning on jumps must fail.

\

On the other hand, by maintaining a finite activity assumption on the jumps
and some conditions on linkage operators, the ellipticity condition has been
relaxed to hypoellipticity by Forster, L\"utkebohmert and Teichmann $[13]$.
Unfortunately, their 'linkage' condition on the jump 
vector fields excludes many examples of financial interest \footnote{Indeed it is easy to see that this condition 
fails in the case where the jumps in the stock are log normal as in the Merton model (or any example in which the L\'evy measure has full support).}. 
The main contribution of $[13]$, in our
view, is to establish new conditions for integrability of the inverse of the
Malliavin covariance matrix $C$ in presence of jumps.  Recent progress in this
direction was also made by Takeuchi $[26]$ who manages to bypass Norris' lemma,
which, in a sense is the bottleneck of the arguments in $[13]$. Thus, noting that
criteria for
 integrability of $C^{-1}$ are available in the literature, and can also be
checked by hand
 in many examples, we show that suitable
integrability of $C^{-1}$ allows to extend a recent result by Arnaudon and
Thalmaier 
$[1]$ and we so obtain non-anticipating Malliavin weights for hypoelliptic
diffusions
with jumps of possibly infinite activity also allowing for the 'linkage' condition in $[13]$ to be relaxed.

    It is worthwhile to ponder for a moment which financial examples really
benefit from BEL / Fournié et al type formulae. The standard hypoelliptic
example in finance is an Asian option but computation of Greeks with
(intelligently chosen) finite difference perform rather well. In fact, most
jump diffusion models used in practice have an essentially\footnote{The gap
between elliptic and
what one has in some real examples is subject of Section 5 of this paper.}
 elliptic diffusion part and also quasi-closed form expressions for European
option prices and the usual Greeks,
typically by Fourier methods i.e. by low-dimensional integration. Thus, the
focus should really be on instruments without (quasi-)closed form prices for
which finite difference methods perform poorly. In fact, there is a very
popular family of such contracts in equity markets, namely digital cliquets,
and the numerical difficulties for risk management are well-known to
practitioners. Surprisingly perhaps, there seems to be no result in the
literature that applies to computing sensitivities of digital cliquets under,
for instance, the Heston-model with jumps: the separability conditions of $[10]$
are far too stringent, the relevant statement in $[13]$, Proposition 1 to be
precise, still contains the (here unnecessary) linkage condition which is not
satisfied \footnote{One could re-run the Malliavin calculus arguments
of $[13]$ in the elliptic
setting to get rid of this condition or, in fact, make rigorous the iterated
conditioning argument outlined above.} nor do Heston-type models
satisfy the strong $C^{\infty}$ assumptions of $[13]$.\footnote{This is just to
say, that approximation arguments
similar to those discussed in Section 5 of this paper would be needed.}

This paper is organised as follows. We prove that the Bismut-Elworthy-Li
formula holds for a generic non-degenerate time-inhomogeneous Markovian
jump-diffusion; $\pi$ is given explicitly as a stochastic integral
involving the flow and the right-inverse of the diffusion matrix, just as
in the classical Bismut-Elworthy-Li formula (which is recovered in the
absence of jumps). A similar presentation is given for the second
derivative. In section $4$ we demonstrate how Malliavin calculus may be used
with an appropriate choice
of perturbation to provide explicit weights in a hypoelliptic setting.
In section 5, as a case study, we show how to represent the spot
sensitivity\footnote{Delta would be a misnomer here. Cont and Tankov $[9]$
contains a nice discussion of how these concepts are related in presence of
jumps.} in the Heston model with jumps (also known as SVJ) and the Matytsin
double jump model (an extension of SVJ and also known as SVJJ). Both models
are described in detail in Gatheral $[16]$ and are popular in the industry
because of their quasi-closed form solutions for European options in terms
of Fourier-transforms, which we use for numerical benchmarks for some
simulations in the last chapter. An honest application of the BEL
formula\footnote{... and Malliavin techniques in general ...}
to these (and many other practically relevant) examples requires approximation
argument which, in our view, have been neglected
in the literature.

The authors would like to thank James Norris, Chris Rogers and Anton Thalmaier for related
discussions.

\section{Preliminaries}

We collect some background material from Gikhman and Skorohod $[17]$. Our
focus is on the strong solution of \begin{equation} \begin{split}
x_{t}^{z} = z + \int_{0}^{t} Z(s, & x_{s-}^{z})ds + \int_{0}^{t} X(s,
x_{s-}^{z})dW_{s}\\ & + \int_{0}^{t} \int_{E} Y(t, x_{s-}^{z},y)
(\mu-\nu)(dy,ds) \end{split} \end{equation} where $W_{t}
\equiv(W_{t}^{1},\ldots, W_{t}^{m})$ is an $\mathbb{R}^{m} $-valued
Brownian motion on some probability space $(\Omega, \mathcal{F}_{t},
\mathbb{P})$ and $\mu$ is a $(\Omega, \mathcal{F}_{t}, \mathbb{P})$-Poisson
random measure on $E\times[0,\infty)$ for some topological space $E$ such
that $\nu$, the compensator of $\mu$, is of the form $G(dy) dt$ for some
$\sigma $-finite measure $G$. The vector fields $Z(t,x)$ and $Y(t,x,y) \in
\mathbb{R}^{d}$, $X(t,x) \in\mathcal{L}(\mathbb{R}^{m},\mathbb{R}^{d})$ for
all $t \in[0,T]$, $x \in\mathbb{R}^{d}$ and $y \in E$. We will always
assume at least the following conditions which guarantee the existence and
uniqueness of a solution to the SDE (see Gikhman and Skorohod $[17]$)

\begin{enumerate}
\item For all $x \in\mathbb{R}^{d}$ and $t \in[0,T]$
\begin{align*}
| Z(t,x)|^{2} + |X(t,x)|^{2} + \int_{E} |Y(t,x,z)|^{2} G(dz) \leq C (1
+|x|^{2})
\end{align*}

\item For all $x,z \in\mathbb{R}^{d}$ and $t \in[0,T]$
\begin{align*}
| Z(t,x)  &  - Z(t,z)|^{2} + |X(t,x) - X(t,z)|^{2}\\
&  + \int_{E} |Y(t,x,y) - Y(t,z,y)|^{2}G(dz) \leq C |x-z|^{2}%
\end{align*}

\end{enumerate}

Throughout we fix the option expiry time $T>0$ and consider a payoff $f :
\mathbb{R}^{d} \rightarrow\mathbb{R}$, and we frequently work with the
process $x_{t}^{t^{\prime},z}$ for $t^{\prime}<t$ defined as the solution
to the SDE 
\begin{equation} \begin{split} x_{t}^{t^{\prime},z} = z + &
\int_{t^{\prime}}^{t} Z(s, x_{s-}^{t^{\prime },z}) ds + \int_{t}^{t}
X(s,x_{s-}^{t^{\prime},z}) dW_{s}\\ & + \int_{0}^{t} \int_{E}
Y(s,x_{s-}^{t^{\prime},z},y) (\mu- \nu)(dy,ds). \end{split} \end{equation}
We will write $f \in C_{b}^{k}(\mathbb{R}^{j})$ to mean that the function
$f: \mathbb{R}^{j} \rightarrow\mathbb{R}$ is $k$-times continuously
differentiable with $f$ and all its derivatives up to order $k$ uniformly
bounded. Our method of proof will rely on ensuring that the function
$u(t,z) : [0,T] \times\mathbb{R}^{d} \rightarrow\mathbb{R}$ given by
$u(t,z) = \mathbb{E} [f(x_{T}^{t,z})]$ satisfies the backward equation of
Kolmogorov. \begin{align*} \mathcal{G}_{t} u(t,z) + \frac{\partial
u}{\partial t}(t,z) = 0. \end{align*} when $f\in
C_{b}^{2}(\mathbb{R}^{d})$. It is well-known that for $f \in
C^{1,2}_{b}([0,T] \times\mathbb{R}^{d})$, $f$ lies in the domain of
$\mathcal{G}$ and \begin{equation} \begin{split} \mathcal{G}_{t} & f =
\sum_{i} Z^{i}(t,x) \frac{\partial f }{\partial x_{i}}(t,x) + \frac{1}{2}
\sum_{i,j} A^{i,j}(t,x) \frac{\partial^{2}f}{\partial x_{i} \partial
x_{j}}(t,x)\\ & + \int_{E} \left( f(t, x + Y(t,x,y)) - f(t,x) - \sum_{i}
Y^{i}(t,x,y) \frac{\partial f}{\partial x_{i}}(t,x) \right) G(dy)
\end{split} \end{equation} where $A = X(t,x)X^{T}(t,x)$ (here $X^{T}$
denotes the transpose of the matrix X). The arguments leading to the proof
of the following theorem may be found in Gikman and Skorohod $[17]$.
\begin{theorem}(Kolmogorov's backward equation) Let $x_{t}^{t',z}$ for
$t'<t\leq T$ represent the solution to the SDE ($3$) with the vector fields
$Z,X$ and $Y$ satisfying the existence and uniqueness conditions. Suppose
further that for every $t \in [0,T]$ and $y \in E$ the following conditions
hold \begin{align} Z(t,.), X(t,.) \text{ and }& Y(t,.,y) \in
C_{b}^{2}(\mathbb{R}^{d}) \end{align} \begin{align} &\text{ The functions }
Z(.,.), \nabla Z(.,.), \nabla^{2} Z(.,.), X(.,.), \nabla X(.,.), \nabla^{2}
X(.,.) \nonumber \\ & \text{ and } \int_{E} |Y(.,.,y)|^{2} G(dy) , \int_{E}
| \nabla Y(.,.,y)|^{2}G(dy) , \int_{E}|\nabla^{2}Y(.,.,y)|^{2} G(dy)
\nonumber \\ & \text{ are all continuous on }[0,T] \times \mathbb{R}^{d}
\end{align} \begin{align} \sup_{(t,x) \in [0,T] \times \mathbb{R}^{d}}
\left( \int_{E} \left( |Y(t,x,y)|^{k} + |\nabla Y(t,x,y)|^{k} \right.
\right. & \big) G(dy) \bigg) < \infty \end{align} for $k=2,3,4$. Then if $f
\in C_{b}^{2}(\mathbb{R}^{d})$ the function $u(t,z) =
\mathbb{E}[f(x_{T}^{t,z})]$ is such that $u \in C_{b}^{1,2}([0,T]\times
\mathbb{R}^{d})$ and satisfies \begin{equation*} \frac{\partial u}{\partial
t}(t,z) + \mathcal{G}_{t} u = 0 \end{equation*} with boundary condition
$\lim_{t \rightarrow T} u(t,z) = f(z)$ and where $\mathcal{G}_{t}$ is given
by $(4)$ . \end{theorem}

\section{The Bismut-Elworthy-Li formula for \newline jump-diffusions}

The argument of Elworthy and Li $[11]$ extends in a straight-forward way to
jump-diffusions. \begin{theorem} Fix some $T>0$ and consider $x_{t}^{z}
\equiv x_{t}^{0,z}$ the solution to SDE $(2)$ on the interval $[0,T]$ and
suppose that the conditions of Theorem $1$ are satisfied. Further assume
that the diffusion matrix $X(t,x) $ has a right inverse $R(t,x)$, and
satisfies the following uniform ellipticity condition \begin{align*} y^{T}
X(t,x) X^{T}(t,x) y \geq \epsilon |y|^{2} \end{align*} for every $t \in
[0,T]$, $x,y \in \mathbb{R}^{d}$ and some $\epsilon > 0$. Then, if $a \in
L^{2}[0,T]$ is any deterministic function which satisfies \begin{align*}
\int_{0}^{T} a(t) dt = 1 \end{align*} and $f\in C_{b}^{2}(\mathbb{R}^{d})$
the following is true for all $1\leq k\leq d$ \begin{equation}
\frac{\partial }{\partial z_{k}} \mathbb{E}[f(x_{T}^{z}) ] = \mathbb{E}
\left[ f(x_{T}^{z}) \int_{0}^{T}a(t) \left( R(t,x^{z}_{t-}) \frac{\partial
x_{t-}^{z}}{\partial z_{k}}\right)^{T}dW_{t} \right]. \end{equation}
Moreover, if we consider $0 < T_{1}\leq \ldots \leq T_{n}\leq T$ and a
function of the form $f(x^{z}_{T_{1}}, \ldots, x^{z}_{T_{n}})$, where $f
\in C_{b}^{2}(\mathbb{R}^{d} \times \ldots \times \mathbb{R}^{d})$ and let
$a \in L^{2}[0,T]$ be a deterministic function satisfying \begin{align*}
\int_{0}^{T_{1}} a(t) dt = 1. \end{align*} Then, for all $1\leq k\leq d$,
the following is true \begin{align} \frac{\partial}{\partial z_{k}}
\mathbb{E}[f( & x^{z}_{T_{1}},\ldots , x^{z}_{T_{n}})] \nonumber \\ & =
\mathbb{E}\left[ f(x^{z}_{T_{1}}, \ldots, x^{z}_{T_{n}})
\int_{0}^{T_{1}}a(t) \left( R(t,x_{t-}^{z}) \frac{\partial
x_{t-}^{z}}{\partial z_{k}}\right)^{T} dW_{t} \right]. \end{align}
\end{theorem}

\begin{subsubaxiom} In the absence of jumps and with $a(t) = T^{-1}$ on
$[0,T]$ we recover the classical Bismut-Elworthy-Li formula.
\end{subsubaxiom} \begin{subsubaxiom} It is easy to see that the uniform
ellipticity condition gives rise to the fact that $\int_{0}^{t}a(s) \left(
R(s,x_{s-}^{z}) \frac{\partial x_{s-}^{z}}{\partial z_{k}}\right)^{T}
dW_{s} $ is a martingale on $[0,T]$. To see this take $t=s$, $x = x_{s-}$
and $ y = R(s,x_{s-})\frac{\partial x_{s-}^{z}}{\partial z_{k}}$ and
observe that \begin{align*} \left|R(s,x_{s-}^{z}) \frac{\partial
x_{s-}^{z}}{\partial z_{k}}\right|^{2} \leq \epsilon^{-1} \left|
\frac{\partial x_{s-}^{z}}{\partial z_{k}}\right|^{2} \,\,\, \, \text{a.s.}
\end{align*} Consequently, \begin{align*} \mathbb{E}\Bigg[
\int_{0}^{T}a(t)^{2}\Big| R(t,x_{t-}^{z}) & \frac{\partial
x_{t-}^{z}}{\partial z_{k}} \Big|^{2} dt \Bigg] \\ & \leq
\epsilon^{-1}\mathbb{E}\left[\sup_{0\leq t\leq T} \left| \frac{\partial
x_{t-}^{z}}{\partial z_{k}}\right|^{2} \right] \int_{0}^{T} a(t)^{2} dt <
\infty. \end{align*} \end{subsubaxiom}

\begin{proof} For $t<T$ we apply It\^o's formula to the function
\begin{align*} u(t,z) = \mathbb{E}[f(x_{T}^{t,z})] = P_{T-t}f(z)
\end{align*} for $t<T$. Since \begin{align*} \mathcal{G}_{t}u +
\frac{\partial u}{\partial t}(t,x) = 0 \end{align*} the ds-term vanishes,
leaving only a constant and the two martingale terms. Letting $t\rightarrow
T$ we find \begin{equation} \begin{split} f(x_{T}^{z}) = & u(0,x) +
\int_{0}^{T}(\nabla u(s,x_{s-}^{z}) X(s,x_{s-}^{z}))^{T} dW_{s} \\ & +
\int_{0}^{T} \int_{E} ( u(s,x_{s-}^{z} + Y(s,x_{s-}^{z},y)) -
u(s,x_{s-}^{z}) ) (\mu-\nu)(ds,dy) \end{split} \end{equation} The integral
featuring above with respect to $\mu-\nu$ is a discontinuous
$L^{2}$-martingale which is orthogonal to the martingale $\int a(s)
\left(R(s,x_{s-}^{z})\frac{\partial x_{s-}^{z}}{\partial z_{k}} \right)^{T}
dW_{s}$. Multiplying by $\int_{0}^{T}a(s) \left( R(s,x_{s-}^{z})
\frac{\partial x_{s-}^{z}}{\partial z_{k}}\right)^{T} dW_{s}$ and using
It\^o's isometry gives the result \begin{align}
\mathbb{E}\left[f(x_{T}^{z}) \int_{0}^{T}a(s) \left( R(s,x_{s-}^{z})
\frac{\partial x_{s-}^{z}}{\partial z_{k}}\right)^{T} dW_{s} \right] & =
\mathbb{E} \left[ \int_{0}^{T}a(s) \nabla u(s,x_{s-}^{z}) \frac{\partial
x_{s-}^{z}}{\partial z_{k}} ds \right] \nonumber \\ & = \int_{0}^{T} a(s)
\mathbb{E} \left[ \nabla u(s,x_{s-}^{z}) \frac{\partial
x_{s-}^{z}}{\partial z_{k}} \right] ds \nonumber \\ & = \int_{0}^{T}a(s)
\frac{\partial}{\partial z_{k}} \mathbb{E} [ u(s,x_{s-}^{z}) ] ds \nonumber
\\ & = \int_{0}^{T}a(s) \frac{\partial}{\partial z_{k}}
\mathbb{E}[f(x_{T}^{z})] ds \nonumber \\ & = \frac{\partial}{\partial
z_{k}} \mathbb{E}[f(x_{T}^{z})] \end{align} We justify the progression from
the second to third line by a routine argument based on the boundedness of
$\nabla u$ and $\nabla^{2}u$ and the definition of $\frac{\partial
x_{t}^{z}}{\partial z_{j}}$ as the $L^{2}$-limit (as $h\rightarrow 0$) of
the random variables $h^{-1}(x_{t}^{z+he_{k}}-x_{t}^{z})$ for fixed $t$.
Also, we justify the third to fourth line in $(11)$ by the observation that
$u(t,x_{t}^{z}) \rightarrow u(s,x_{s-})$ almost surely as $t\uparrow s$ and
so bounded convergence gives $\mathbb{E}[u(t,x_{t}^{z})] \rightarrow
\mathbb{E}[u(s,x_{s-})]$. But for each $t\in[0,T]$we have
$\mathbb{E}[u(t,x_{t}^{z})] = \mathbb{E}[\mathbb{E}[f(x_{T}^{z}) |
\mathcal{F}_{t}]] = \mathbb{E}[f(x_{T}^{z})]$, so $\mathbb{E}[u(s,x_{s-})]
= \mathbb{E}[f(x_{T}^{z})]$. For the final part, we note that the function
$g: \mathbb{R}^{d} \rightarrow \mathbb{R}$ defined by $g(x) =
\mathbb{E}[f(x,x_{T_{2}}^{T_{1},x},\ldots x_{T_{n}}^{T_{1},x})]$ has the
property that $g \in C_{b}^{2}(\mathbb{R}^{d})$ and, moreover, by the
Markov property \begin{align*} g(x_{T_{1}}^{z}) & =
\mathbb{E}[f(x_{T_{1}}^{z},\ldots, x_{T_{n}}^{z}) | \sigma(x_{T_{1}}^{z}) ]
= \mathbb{E}[f(x_{T_{1}}^{z},\ldots, x_{T_{n}}^{z}) | \mathcal{F}_{T_{1}}]
\,\,\, \text{a.s.} \end{align*} Consequently, by $(8)$ we have
\begin{align*} \frac{\partial}{\partial z_{k}} \mathbb{E}
[f(x^{z}_{T_{1}},\ldots , x^{z}_{T_{n}})] & = \frac{\partial}{\partial
z_{k}} \mathbb{E}[g(x_{T_{1}}^{z})] \\ & = \mathbb{E} \left[
g(x_{T_{1}}^{z}) \int_{0}^{T_{1}}a(t) \left( R(t,x^{z}_{t-}) \frac{\partial
x_{t-}^{z}}{\partial z_{k}}\right)^{T}dW_{t} \right], \end{align*} which
concludes the proof. \end{proof}

\begin{subsubaxiom} Under stronger condition on the vector fields ( see
Theorem ($2$-$28$) of Bichteler, Jacod and Gravereaux $[5]$ ) we can ensure
the existence of a density $p_{T}(z,y)$ for the random variable $x_{T}^{z}$
with $p_{T} \in C^{1}(\mathbb{R}^{d} \times \mathbb{R}^{d})$. It is then
possible to relax the regularity restrictions on $f$ so that we need only
make a measurability assumption on $f$. \end{subsubaxiom}

\begin{subsubaxiom} The result $(9)$ extends the result of Section $3.2$
Fourni\'e et al $[14]$ to jump-diffusions. We notice that the form of the
results do not correspond exactly, their weight is represented by
\begin{align*} \pi = \int_{0}^{T} a(t) \left( R(t,x^{z}_{t-})
\frac{\partial x_{t-}^{z}}{\partial z_{k}} \right)^{T} dW_{t}, \end{align*}
with $a \in L^{2}[0,T]$ satisfying $\int_{0}^{T_{i}} a(t) =1$ for all
$1\leq i\leq n$, and our weight $\tilde{\pi}$ is a particular case of this
when $a =0$ on $[T_{1},T]$. However, it is clear that if $a \not= 0$ on
$[T_{1},T]$ then \begin{align*} Var(\tilde{\pi}) \leq Var(\pi).
\end{align*} Since the efficiency of Monte Carlo is optimised by the choice
of the minimal variance weight we would always choose $a\equiv 0$ on
$[T_{1},T] $ and hence there is no conceivable practical advantage to
representing the weight by $\pi$. \end{subsubaxiom} We may adapt this
approach to deal with higher order derivatives as well.

\begin{theorem} Suppose that $XX^{T}$ is uniformly elliptic and further
assume that the conditions on the vector fields are strengthened so that
the following conditions are satisfied. For every $t \in [0,T]$ and $y \in
E$ \begin{align} Z(t,.), X(t,.) \text{ and } Y(t,.,y) \in C_{b}^{\infty}
(\mathbb{R}^{d}). \end{align} \begin{align} & \text{For every $l \in
\mathbb{N} \cup {\{0\}}$, } \nabla^{l} Z(.,.), \nabla^{l}X(.,.) \text{ and
} \int_{E} |\nabla^{l} Y(.,.,y)|^{2} G(dy) \nonumber \\ & \text{are
continuous on} [0,T] \times \mathbb{R}^{d} \end{align} \begin{align} \text{
For}& \text{ $r \in \mathbb{N}$ with $r\geq 2$ and $l=1,2$.} \nonumber \\ &
\sup_{(t,x) \in [0,T] \times \mathbb{R}^{d}} \left( \int_{E}(
| Y(t,x,y)|^{r} + |\nabla^{l} Y(t,x,y)|^{r}) G(dy) \right) < \infty
\end{align} Then, if $f \in C_{b}^{3}(\mathbb{R}^{d})$ and, for each $t \in
[0,T]$, $R(t,.) \in C^{1}_{b}(\mathbb{R}^{d})$ (where the bounds on
$R(t,.)$ and $\nabla R(t,.)$ hold uniformly in $t \in [0,T]$, the following
formula holds for all $1\leq j,k\leq d$ \begin{align*}
\frac{\partial^{2}}{\partial z_{j} \partial z_{k}} \mathbb{E}[f(x_{T}^{z})]
&\\ = \frac{4}{T^{2}} \mathbb{E}\Bigg[ f(x_{T}^{z}) \int_{T/2}^{T} &
\left(R(t,x_{t-}^{z}) \frac{\partial x_{t-}^{z}}{\partial z_{j}}\right)^{T}
dW_{t} \int_{0}^{T/2} \left(R(t,x_{t-}^{z}) \frac{\partial
x_{t-}^{z}}{\partial z_{k}}\right)^{T} dW_{t}\Bigg] \\ & +
\frac{2}{T}\mathbb{E}\left[ f(x_{T}^{z})\int_{0}^{T/2} \left(\nabla
R(t,x_{t-}^{z}) \frac{\partial x_{t-}^{z}}{\partial z_{j}} \frac{\partial
x_{t-}^{z}}{\partial z_{k}} \right)^{T} dW_{t}\right] \\ & + \frac{2}{T}
\mathbb{E}\left[f(x_{T}^{z}) \int_{0}^{T/2}\left(
R(t,x_{t-}^{z})\frac{\partial^{2} x_{t-}^{z}}{\partial z_{j}\partial z_{k}}
\right)^{T}dW_{t} \right]. \end{align*} Moreover, if we consider $0 <
T_{1}\leq \ldots \leq T_{n}\leq T$ and a function of the form
$f(x^{z}_{T_{1}}, \ldots, x^{z}_{T_{n}})$, where $f \in
C_{b}^{3}(\mathbb{R}^{d} \times \ldots \times \mathbb{R}^{d})$. Then, the
above result remains true when we replace $f(x^{z}_{T})$ by
$f(x^{z}_{T_{1}}, \ldots, x^{z}_{T_{n}})$ and $T$ by $T_{1}$ in the above
formula. \end{theorem} \begin{subsubaxiom} The conditions on the vector
fields are stronger than needed, but we state them in their current form
for simplicity. \end{subsubaxiom} \begin{proof} Define the function $w :
\mathbb{R}^{d}\times \mathbb{R}^{d} \rightarrow \mathbb{R}$ by $w(x,y) =
\nabla f(x)y$. Then it is easy to verify  that the
function $p : [0,T] \times \mathbb{R}^{d} \times \mathbb{R}^{d} \rightarrow
\mathbb{R}$ defined by $p(t,z,y) = \mathbb{E} \left[ w\left(x_{T}^{t,z},
\frac{\partial x_{T}^{t,x}}{\partial x_{j}}\right)\right]$ satisfies the
backward equation associated to the generator of the $\mathbb{R}^{d} \times
\mathbb{R}^{d}$-valued diffusion $\left(x_{t}^{z},\frac{\partial
x_{t}^{z}}{\partial z_{j}}\right)^{T}$. The argument now proceeds as
before; applying It\^o's formula to $p\left(t,x_{t}^{z},\frac{\partial
x_{t}^{z}}{\partial z_{j}} \right)$, letting $t\rightarrow T$ and then
multiplying by $\int_{0}^{T}\left( R(t,x_{t-}^{z}) \frac{\partial
x_{t}^{z}}{\partial z_{k}}\right)^{T} dW_{t}$ and taking expectations
allows the argument to be concluded as in Theorem $2.3$ of Elworthy and Li
$[11]$. \end{proof}

\section{Relaxing the Ellipticity Criterion}
For simplicity we now assume that the vector fields are time homogeneous.  
We denote by $U$ the $\mathcal{L}(\mathbb{R}^{d},\mathbb{R}^{d})$-valued 
process given by $U_{t} = \nabla_{z} x_{t}^{z}$ and denote its inverse, when it exists,  by $V_{t}$.
We define the Malliavin covariance matrix 
\begin{align*}
C_{t}(z) = \int_{0}^{t} (V_{s-}X(x_{s-}^{z}))(V_{s-}X(x_{s-}^{z}))^{T} ds
\end{align*}
and make the following a standing assumption.
\begin{subsubsubtheorem}
For fixed $T>0$, $C_{T}$ is invertible a.s. and moreover $|C_{T}^{-1}| \in L^{p}$ for all $p\geq 1$. 
\end{subsubsubtheorem}
This assumption is known to be true in certain cases, for instance it holds in the diffusion case
under H\"ormander conditions on the vector fields (see Nualart $[21]$), and more recently it has been shown to hold 
in the jump diffusion case for finite intensity jumps under uniform H\"ormander condition (see Forster, L\"utkebohmert and Teichmann $[13]$).
For more general jump processes the problem is more involved but ideas in this setting have been developed  
in Cass $[8]$ and Takeuchi $[26]$.

We now prove an extension of Theorem $3.2$ of Arnaudon and Thalmaier $[1]$ which allows us to give an explicit representation of the weight in terms
of an adapted $\mathbb{R}^{d}$-valued process .
Note that the result of Forster, L\"utkebohmert and Teichmann $[13]$ where the weight is given in the form of a anticipating 
Skorokhod integral may be converted into sum of integrals of adapted processes using the expansion formula $(1.49)$ in Nualart $[21]$.
A representation of this type is more desirable from the point of view of simulation.
First we recall some concepts from Malliavin calculus. Let $a$ be an $\mathcal{L}(\mathbb{R}^{d},\mathbb{R}^{m})$-valued
previsible process such that for $T>0$ fixed
\begin{align}
\mathbb{E}\left[ \exp\left( \frac{1}{2}\int_{0}^{t} |a_{s} h|^{2} ds \right) \right] < \infty , \,\, \text{$h \in \mathbb{R}^{d}$ locally at $0$}
\end{align}
and define a new probability measure by 
\begin{align*}
Z_{T}^{h} = \frac{d \mathbb{P}^{h}}{d \mathbb{P}} \Big{|}_{\mathcal{F}_{T}} = \exp\left( - \int_{0}^{t}< a_{s}h,dW_{s}>
 -\frac{1}{2}\int_{0}^{t} |a_{s} h|^{2}ds \right),
\end{align*}
 and let $Z_{t} = \mathbb{E}[Z_{T}^{h} | \mathcal{F}_{t}]$ for $0\leq t\leq T$.
Introduce a perturbed Brownian motion $dW_{t}^{h} = dW_{t}+ a_{s}h dt$ and denote by  $x_{t}^{h}$,$C_{t}^{h}(z)$  the 
processes corresponding to $x_{t}^{z}$ and $C_{t}(z)$ when the underlying Brownian motion is replaced by $W_{t}^{h}$. The crucial ingredient 
to the following result is the observation that the change of measure decribed above has no effect on the Poisson random measure.
\begin{theorem}
Suppose Assumption $1$ is in force along with the following conditions on the vector fields
\begin{align*}
& Z(.) \in C_{b}^{\infty}(\mathbb{R}^{d},\mathbb{R}^{d}), X(.) \in C_{b}^{\infty}(\mathbb{R}^{d},\mathcal{L}(\mathbb{R}^{m},\mathbb{R}^{d})),  
Y(.,y) \in C_{b}^{\infty}(\mathbb{R}^{d},\mathbb{R}^{d}) \\
&\text{and} \, \sup_{x \in \mathbb{R}^{d}} \int_{E} |\nabla^{n}_{x}Y(x,y)|^{2}G(dy) < \infty \, \, \text{for all} \, n \in \mathbb{N}.
\end{align*}
Further assume 
\begin{align}
\sup_{x \in \mathbb{R}^{d}}\sup_{y \in E} | (I+ \nabla_{x}Y(x,y))^{-1}| < \infty.
\end{align}
Then, for any $f \in C^{1}_{c}(\mathbb{R}^{d})$ and $j \in \{1,2,\ldots,d\}$ we have
\begin{align*}
 \frac{\partial}{\partial z_{j}} \mathbb{E}[f(x_{T}^{z})]  
= \mathbb{E}\Bigg[ & f(x_{T}^{z})\Bigg(\Big(\int_{0}^{T}V_{t-}X(x_{t-})dW_{t}
\Big)^{T}C_{T}^{-1}(z)e_{j}  \\
&+ \sum_{k=1}^{d} \Bigg( C_{T}^{-1}(z) \Big( \frac{\partial}{\partial h_{k}}
\Big{|}_{h=0} C_{T}^{h}(z) \Big) C_{T}^{-1}(z) \Bigg)_{k,j}\Bigg) \Bigg]
\end{align*}
\end{theorem}
\begin{subsubaxiom}
Condition $(16)$ is there to ensure both the existence of $V = U^{-1}$ and that $V \in L^{p}$ for all $p\geq 1$. In practice this can often be
relaxed in favour of some less stringent condition (see Example $2$ below).
\end{subsubaxiom}
\begin{proof}
The fact that $\mu$ remains a Poisson random measure with compensator $\nu$ under $\mathbb{P}^{h}$ follows from Theorems $(3.15)$ and
$(3.34)$
of Jacod $[20]$ . We then observe,since $x$ is a strong solution to the 
SDE $(2)$, that
\begin{align}
\sum_{k=1}^{d} \frac{\partial}{\partial h_{k}} \Big{|}_{h=0} \mathbb{E}[f(x_{T}^{h}) Z_{T}^{h} (C_{T}^{h}(z)^{-1})_{k,j}] = 0. 
\end{align}
Choosing the perturbation 
\begin{align*}
a_{s-}^{n} =  V_{s-}X(x_{s-})1_{\{s\leq \tau_{n}\}}
\end{align*}
with an increasing sequence of previsible stopping times ($\tau_{n}$) chosen such that $a^{n}$ satisfies 
condition $(15)$ and such that $\tau_{n}\uparrow T$.An elementary application of It\^o's formula can be used to show 
\begin{align*}
\frac{\partial}{\partial h_{k}} \Big{|}_{h=0} x_{T}^{h} = U_{T} \int_{0}^{\tau_{n}}V_{s-}X(x_{s-})(V_{s-}X(x_{s-}))^{T}ds e_{k} =U_{T}C_{\tau_{n}}e_{k}.
\end{align*}
Using this we may expand $(16 )$ to get 
\begin{align*}
\mathbb{E}[f(x_{T}^{z})U_{T}C_{\tau_{n}}C_{T}^{-1}e_{j}] = \mathbb{E}\Bigg[ & f(x_{T}^{z})\Bigg(\Big(\int_{0}^{\tau_{n}}
V_{t-}X(x_{t-})dW_{t}
\Big)^{T}C_{T}^{-1}(z)e_{j}  \\
&- \Big( \sum_{k=1}^{d}  \frac{\partial}{\partial h_{k}}
\Big{|}_{h=0}(C_{T}^{h})^{-1}e_{k}  \Big)^{T} e_{j}\Bigg) \Bigg].
\end{align*}
We let $n \rightarrow \infty$  and expand the second term on the right hand side to give the stated result.
\end{proof}
\begin{subsubsubaxiom}(Bachelier with jumps,  Asian options)
We assume 
\begin{align*}
dS_{t} & =  \sigma dW_{t} + dN_{t} \\
dA_{t} & = S_{t-}dt,
\end{align*} 
with some Poisson process $N$ of finite rate. The Malliavin covariance matrix has the particularly simple form 
\begin{align*}
C_{T} = \sigma^{2} \begin{pmatrix} T & -T^{2}/2 \\ -T^{2}/2 & T^{3}/3 \end{pmatrix} 
\end{align*}
and so the second term on the right hand side of the formula in Theorem $4$ drops out leaving us with
\begin{align*}
\frac{\partial}{\partial S_{0}} \mathbb{E}[f(S_{T},A_{T})] 
= \frac{6}{\sigma T} \mathbb{E}\left[ \left( \frac{1}{T} \int_{0}^{T} W_{t} dt - \frac{1}{3} W_{T}\right)f(S_{T},A_{T}) \right].
\end{align*}
\end{subsubsubaxiom}
\begin{subsubsubaxiom}(Exponential L\'evy,  Asian options) 
Consider the following model for the evolution of a stock price 
\begin{align*}
dS_{t} & = \beta  S_{t-}dt + \sigma S_{t-}dW_{t} + \int_{y \geq  -1 } y S_{t-} (\mu-\nu)(dy,dt) \\
dA_{t} & = S_{t-} dt,
\end{align*} 
with $A_{0} = 0$, and where $W$ is a Brownian motion and $\mu$ a Poisson random measure with compensator $\nu(dy,dt) =G(dy)dt$. 
Make the assumptions that for all $p\geq 1$ and arbitrary $\delta >0$
\begin{align}
\int_{-1}^{-\delta} (1+y)^{-p} G(dy) < \infty \,\,\, \text{ and } \int_{y\geq 1} (1+y)^{p} G(dy) < \infty . 
\end{align}
Condition $(16)$ is not satisified in this example, however it is easy to show by truncating the jumps at some arbitrary level
that the theorem may be applied. Assumption $(18)$ may then be invoked to guarantee the resulting formula remains valid in the limit
as the truncation parameter goes to zero.
We notice also in this case that the vector fields are not bounded and similar approximation results are needed, details
on how this type of argument can be made rigorous are given in the next section but we omit them here for the purpose of clear exposition.  
The Malliavin covariance matrix may be computed 
\begin{align*}
C_{T} = \int_{0}^{T} \begin{pmatrix} \sigma^{2} S_{0}^{2} &  -\sigma^{2} S_{0} A_{t}  \\  -\sigma^{2} S_{0} A_{t}   &  \sigma^{2} A_{t}^{2} \end{pmatrix}dt 
= \begin{pmatrix} \sigma^{2}  S_{0}^{2}T &  -\sigma^{2} S_{0} \int_{0}^{T} A_{t}dt  \\  -\sigma^{2} S_{0} \int_{0}^{T} A_{t}dt  & \sigma^{2} \int_{0}^{T} A_{t}^{2} dt \end{pmatrix}.
\end{align*}
It is easy to show that
\begin{align*}
\frac{\partial}{\partial h_{1}}\Big{|}_{h=0} S_{t}^{h} = \sigma^{2} S_{0} S_{t}t \,\, \,\, \,\,
\frac{\partial}{\partial h_{2}} \Big{|}_{h=0} S_{t}^{h} = -\sigma^{2}S_{t} \int_{0}^{t} A_{s} ds 
\end{align*}
and then
\begin{align*}
\frac{\partial}{\partial h_{i}}\Big{|}_{h=0} A_{t}^{h} = \begin{cases} \sigma^{2}S_{0} \int_{0}^{t} sS_{s-}ds  & \text{ if $i=1$} \\
\sigma^{2} \left( \int_{0}^{t} A_{s}^{2} ds - A_{t}\int_{0}^{t}A_{s} ds  \right) & \text{if $i=2$}\end{cases}.
\end{align*}
We notice that det$C_{t} = \sigma^{4}S_{0}^{2} \left( t \int_{0}^{t} A_{t}^{2} dt - \left(\int_{0}^{t} A_{t} dt \right)^{2} \right)$, and 
\begin{align*}
C_{T}^{-1} =(\text{det}C_{T})^{-1}  
\begin{pmatrix}\sigma^{2} \int_{0}^{T} A_{t}^{2} dt  &  \sigma^{2} S_{0} \int_{0}^{T}A_{t}dt  
\\ \sigma^{2} S_{0} \int_{0}^{T}A_{t}dt  & \sigma^{2}  S_{0}^{2}T \end{pmatrix}.
\end{align*}
We must show that $C_{T}^{-1} \in L^{p}$ for all $p\geq 1$. To see this it suffices to check that 
\begin{align*}
\mathbb{P}\left(\text{det} C_{T} \leq \epsilon \right) \, \, \, \, \text{ is } \,\, o(\epsilon ^{p}) \,\,\, \text{as} \,\, \epsilon \rightarrow 0.
\end{align*}  
To this end we note that for any $ 0< \delta < 1$
\begin{align*}
\mathbb{P}\left(\text{det} C_{T} \leq \epsilon \right) \leq \mathbb{P}& \left(\text{det} C_{T} \leq \epsilon, \inf_{0\leq t\leq T} S_{t} > \delta,
\sup_{0\leq t\leq T} S_{t} < \delta^{-1} \right) \\
& + \mathbb{P} \left( \inf_{0\leq t\leq T} S_{t} \leq \delta \right) + \mathbb{P}\left( \sup_{0\leq t\leq T} S_{t} \geq  \delta^{-1} \right).
\end{align*}
We now show that on the set $ B : =\{\inf_{0\leq t\leq T} S_{t} > \delta,
\sup_{0\leq t\leq T} S_{t} < \delta^{-1} \}$ we have, for $\delta = \delta(\epsilon)$ appropriately chosen , det$C_{T}>\epsilon$. To see this note that
on $B$ we have $\delta t \leq A_{t}(\omega) \leq \delta^{-1} t$ for all $t \in [0,T]$, and so we define
\begin{align*}
\mathcal{A} = \{ f : [0,T] \rightarrow \mathbb{R}, \,\, \text{such that} \,\, \delta t \leq f(t) \leq \delta^{-1}t \,\, \text{for all} \, t \in[0,T] \}.
\end{align*}
Then by examining the form of the determinant we have 
\begin{align*}
det C_{T} \geq \sigma^{4} S_{0}^{2} T^{2} \inf_{f \in \mathcal{A}} \text{var}f(U) 
\end{align*}
where $U\sim \text{Uniform}[0,T]$. We may bound the left hand side from below by Chebyshev's inequality, so that for any $a>0$
\begin{align*}
\text{var}f(U) & \geq a^{2} \mathbb{P}(|f(U) - \mathbb{E}[f(U)]| \geq a) \\
& \geq a^{2} \mathbb{P}(f(U) \leq -a + \mathbb{E}[f(U)] ) \\
& \geq a^{2} \mathbb{P}\left(  f(U) \leq -a + \frac{\delta T}{2} \right) 
\end{align*}
and taking $a = \frac{\delta T}{4}$ gives var$f(U) \geq \frac{\delta^{2} T^{2}}{16} \mathbb{P}\left( f(U) \leq \frac{\delta T}{4} \right)$. Since
$f(t) \leq \delta^{-1} t$ in $[0,T]$ we have $\mathbb{P}\left(f(U) \leq \frac{\delta T}{4}\right) \geq  \mathbb{P}\left( U \leq \frac{\delta^{2}T}{4} \right)$
which gives 
\begin{align*}
\text{var} f(U) \geq \frac{\delta^{4}T^{3}}{64} 
\end{align*}
and so det$C_{T} \geq \frac{\delta^{4}T^{5}\sigma^{4}S_{0}^{2}}{64} : = C\delta^{4}$. Choosing $\delta = (C^{-1}\epsilon)^{1/4}$ we see
that det$C_{T}\geq \epsilon$ on $B$. It therefore suffices to show that 
\begin{align*}
\mathbb{P}\left( \inf_{0\leq t\leq T} S_{t} \leq \epsilon \right) \,\, \text{and}\,\, \mathbb{P}\left( \sup_{0\leq t\leq T} S_{t} > \epsilon^{-1} \right)  
\,\, \text{ are } \,\, o(\epsilon ^{p}) \,\,\, \text{as} \,\, \epsilon \rightarrow 0
\end{align*}
for every $p\geq 1$. We show this for the infimum, the supremum being a simple modification of this argument.
To this end 
we write $S_{t} = S_{0}e^{X_{t}}$ where $X_{t}$ is a L\'evy process with triplet $(\sigma^{2},\tilde{\mu},\tilde{G})$
with 
\begin{align*}
\tilde{\mu}  = \mu - \frac{1}{2} \sigma^{2} - \int_{|y| \geq 1}y G(dy) , \,\,
\tilde{G}(A)  = G( \{ e^{x}-1 : x \in A\}) \,\,\, \text{for $A \in \mathcal{B}(A)$}.
\end{align*}
It is easy to verify using the definition of $\tilde{G}$ that and assumptions $(18)$ that
$\int_{|x|\geq 1} e^{px} \tilde{G}(dx) < \infty$ for all 
$p \in \mathbb{R}$ and from Theorem $25.17$ of Sato $[25]$ this means that $\mathbb{E}[e^{pX_{t}}] < \infty$ and, moreover,
$\mathbb{E}[e^{pX_{t}}] = e^{t \Psi (p)}$ where 
\begin{align*}
\Psi (u) = \frac{1}{2}\sigma^{2}u^{2} + \int_{\mathbb{R}} (e^{ux} - 1 - ux 1_{[-1,1]}(x)) \tilde{G}(dx) + \tilde{\mu}u.
\end{align*}
Next, by choosing $a> 0$ sufficiently large we may ensure that
$\mathbb{P} \left( \sup_{0\leq t\leq T} X_{t} > a \right) \leq 1/2 $. Then, 
\begin{align*}
\mathbb{P}(\inf_{0\leq t\leq T} X_{t} \leq -2a )   = \mathbb{P}( \inf_{0\leq t\leq T} X_{t}& \leq -2a , X_{T} \leq -a)  \\
 & + \mathbb{P}(\inf_{0\leq t\leq T} X_{t} \leq -2a , X_{T} > -a) 
\\   \leq \mathbb{P}(X_{T} \leq -a&)   + \mathbb{P}(\inf_{0\leq t\leq T} X_{t} \leq -2a , X_{T} > -a),  
\end{align*}
and for the second term in the preceeding inequality we may use the strong Markov property at the stopping time
$\, \zeta = \inf\{ t \geq 0 : X_{t} \leq -2a\}$ to give  $\, \mathbb{P}(\inf_{0\leq t\leq T} X_{t} \leq -2a , X_{T} > -a) \leq 1/2 \mathbb{P}(
\inf_{0\leq t\leq T} X_{t} \leq -2a)$, and so for any $p\geq 1$ we have
\begin{align*}
\mathbb{P}(\inf_{0\leq t\leq T} X_{t} \leq -2a ) \leq 2\mathbb{P}(X_{T}\leq -a) \leq 2e^{-pa} e^{T\Psi (-p)}.
\end{align*}
Finally we finish by noting that for $\epsilon$ sufficiently small 
\begin{align*}
  \mathbb{P} ( \inf_{0\leq t\leq T} S_{t} \leq \epsilon) & = \mathbb{P}\left( \inf_{0\leq t \leq T}X_{t} \leq \log \left( \frac{\epsilon}{S_{0}}\right) \right) 
\leq \frac{2 e^{T\Psi (-p)}}{S_{0}^{p}} \epsilon^{p}.
\end{align*}
Using these facts and the previous theorem we have a random variable  $\pi = \pi_{1} + \pi_{2}$, where 
\begin{align*}
\pi_{1} = \frac{1}{\sigma S_{0}}\left( \frac{W_{T} \int_{0}^{T} A_{t}^{2}dt  - 
\int_{0}^{T} A_{t}dW_{t}\int_{0}^{T}A_{t} dt}{T\int_{0}^{T} A_{t}^{2} dt - \left(\int_{0}^{T} A_{t} dt \right)^{2}} \right)
\end{align*}
and $\pi_{2} = \pi_{2,1} + \pi_{2,2}$ with 
\begin{align*}
\pi_{2,1} = \frac{S_{0}^{2}\sigma^{6}}{(\text{det}C_{T})^{2}}  \Bigg( -2 &
\int_{0}^{T} A_{t}^{2}dt \int_{0}^{T} A_{t} dt \int_{0}^{T} \frac{\partial}{\partial h_{1}}\Big{|}_{h=0} A_{t}^{h} dt \\
& + \left(\int_{0}^{T} A_{t} dt\right)^{2} \int_{0}^{T} \frac{\partial}{\partial h_{1}}\Big{|}_{h=0} (A_{t}^{h})^{2} dt \Bigg)
\end{align*} 
and 
\begin{align*}
\pi_{2,2}& = \frac{S_{0}^{3} \sigma^{6}}{(\text{det}C_{T})^{2} }
\Bigg( - \left(\int_{0}^{T} A_{t} dt \right)^{2} \int_{0}^{T} \frac{\partial}{\partial h_{2}} \Big{|}_{h=0} A_{t}^{h} dt \\
& + T\left( \int_{0}^{T} A_{t} dt \int_{0}^{T} \frac{\partial}{\partial h_{2}} \Big{|}_{h=0} (A_{t}^{h})^{2} dt 
-  \int_{0}^{T} A_{t}^{2} dt \int_{0}^{T} \frac{\partial}{\partial h_{2}} \Big{|}_{h=0} A_{t}^{h} dt \right) \Bigg)
\end{align*}
such that 
\begin{align*}
\frac{\partial}{\partial S_{0}} \mathbb{E}[f(S_{T},A_{T})] = \mathbb{E}[f(S_{T},A_{T}) \pi].
\end{align*}
Numerical implementation of these results shows a good degree of accuracy comparable to that achieved by finite difference Monte Carlo in the
case of a European call. 
\end{subsubsubaxiom}

\section{Examples}

We show how the formula derived in the previous section should be
implemented to obtain appropriate representations. We will find that the
restrictions imposed by Theorems $2$ and $3$ on the vector fields are often
too stringent and that we have to get round this problem by localisation.

\subsection{Stochastic volatility models with jumps}

We will consider a volatility process $\sigma_{t}$ described by the Heston
model \begin{align} d \sigma_{t}^{2} = \kappa(\theta- \sigma_{t}^{2}) dt +
\eta\sigma_{t} dW_{t}. \end{align} We will need the following lemma
\begin{axiom} For any parameter choice with $2 \kappa \theta > \eta^{2}$
and for every finite $T>0$ \begin{align*} \sup_{0\leq t\leq T}
\mathbb{E}[\sigma_{t}^{-2}] < \infty. \end{align*} \end{axiom}
\begin{proof} We let $Y_{t}$ be the squared $\delta$-dimensional Bessel
process defined as the unique strong solution to the SDE \begin{align*}
Y_{t} = \sigma_{0}^{2} + \delta t + 2 \int_{0}^{t} \sqrt{Y_{s}} dW_{s}.
\end{align*} For the choice $\delta = \frac{4 \kappa \theta}{\eta^{2}}$ we
can relate $Y_{t}$ and $\sigma_{t}^{2}$ by the time change (see
Going-Jaeschke and Yor $[19]$) \begin{align*} \sigma_{t}^{2} = e^{-\kappa t}
Y \left( \frac{\eta^{2}}{4 \kappa} (e^{\kappa t} -1) \right). \end{align*}
If we let $f(t) = \frac{\eta^{2}}{4 \kappa} (e^{\kappa t} -1)$ and $T^{*} =
f(T) < \infty$ then as $t$ takes values in $[0,T]$ so $f(t)$ ranges over
$[0,T^{*}]$. Consequently, \begin{align} \sup_{0\leq t\leq T}
\mathbb{E}[\sigma_{t}^{-2}] = \sup_{0\leq t\leq T} \mathbb{E}[(e^{-\kappa
t}Y(f(t)))^{-1}] \leq e^{\kappa T} \sup_{0\leq t\leq T^{*}} \mathbb{E}[
Y_{t}^{-1}]. \end{align} Next we notice from the expression for the Laplace
transform of $Y_{t}$ (Revuz and Yor $[23]$, page 422) \begin{align*}
\mathbb{E}[ Y_{t}^{-1}] = \int_{0}^{\infty} \mathbb{E}[e^{-\lambda Y_{t}}]
d \lambda = \int_{0}^{\infty} (1+ 2\lambda t)^{-\delta/2} \exp \left(
\frac{-\lambda \sigma_{0}^{2}}{1+2\lambda t} \right) d \lambda.
\end{align*} From this and the fact that $\delta > 2$ we have, for any
$\epsilon > 0$, \\ $\sup_{\epsilon\leq t\leq T^{*}} \mathbb{E}[Y_{t}^{-1}]
\leq \int_{0}^{\infty} (1+ 2\lambda t)^{- \delta / 2} d\lambda \leq
(2\epsilon)^{-1}$. So the proof will be complete if we can show
\begin{align*} \limsup_{t\rightarrow 0} \int_{0}^{\infty} (1+ 2\lambda
t)^{-\delta/2} \exp \left( \frac{-\lambda \sigma_{0}^{2}}{1+2\lambda t}
\right) d \lambda < \infty. \end{align*} 
 By using the substitution $\xi = 1- (1+2 \lambda t)^{-1}$
and writing $z = t^{-1}$, $y = \frac{\sigma_{0}^{2}z}{2}$  we need to examine the behaviour of
\begin{align*}
 \frac{z}{2} \int_{0}^{1} (1-\xi)^{\delta /2 -2}
e^{-\sigma_{0}^{2}z \xi /2} d \xi & = \sigma_{0}^{-2} 
 \int_{0}^{1} (1-\xi)^{\delta /2 -2} ye^{- y \xi } d \xi  \end{align*} 
as $y \rightarrow \infty$, and it suffices the check that the expression on the right hand side is bounded 
for large $y$.  To show this, first suppose $\frac{\delta}{2} -2 \geq  0$ then we trivially have
\begin{align*}
\sigma_{0}^{-2} 
 \int_{0}^{1} (1-\xi)^{\delta /2 -2} ye^{- y \xi } d \xi \leq \sigma_{0}^{-2}(1-e^{-y}) \leq \sigma_{0}^{-2}.
\end{align*}
Next, suppose $\frac{\delta}{2} - 2 < 0$, then by making the substitution  
 $w = (1-\xi)y$ and noticing that for $y>1$
\begin{align*}
\int_{0}^{1} (1-\xi)^{\delta /2 -2} ye^{- y \xi } d \xi  
& =  e^{-y} \int_{0}^{y} w^{\delta/2 -2} e^{w} dw \\
& \leq e^{-y}\left( e^{1} \int_{0}^{1} w^{\delta/2 - 2} dw + \int_{1}^{y} e^{w} dw \right) ,
\end{align*}
we see that  the right hand side may be bounded uniformly in $y$ since $\frac{\delta}{2} -2 > -1$.
\end{proof}

It will be convenient to think of the process $\sigma_{t}$ instead, so
writing $X_{t} = \log S_{t}$ to represent the evolution of the logarithm of
the stock price the system can be described by the vector SDE \[
\begin{split} \begin{pmatrix} X_{t}\\ \sigma_{t} \end{pmatrix} =
\begin{pmatrix} x\\ \sigma_{0} \end{pmatrix} & + \int_{0}^{t}
\begin{pmatrix} r - \frac{1}{2}\sigma_{s}^{2}\\ \left(
\frac{\kappa\theta}{2} - \frac{\eta^{2}}{8} \right) \frac{1} {\sigma_{s}}
- \frac{\kappa}{2} \sigma_{s} \end{pmatrix} ds + \int_{0}^{t}
\begin{pmatrix} \sqrt{1-\rho^{2}} \sigma_{s}\\ 0 \end{pmatrix} dZ_{s}\\ & +
\int_{0}^{t} \begin{pmatrix} \rho\sigma_{s}\\ \frac{\eta}{2} \end{pmatrix}
dW_{s} + \int_{0}^{t} \int_{E} \begin{pmatrix} y\\ 0 \end{pmatrix} (\mu-
\nu) (dy,ds). \end{split} \] We shall call this model SVJ. Before the next
theorem we introduce the notation $C_{c}^{k}(\mathbb{R}^{d})$ to indicate
the set of real-valued, $k$-times differentiable , compactly supported
functions with domain $\mathbb{R}^{d}$. We then define $I(\mathbb{R}^{d})$
to be the collection of indicator functions of the form $1_{(a,b)},
1_{(a,b]}, 1_{[a,b)}$ or $1_{[a,b]}$ for some $|a| < |b| < \infty$ and
finally a class of real-valued functions on $\mathbb{R}^{d}$,
$\mathcal{J}(\mathbb{R}^{d})$, by \begin{align*}
\mathcal{J}(\mathbb{R}^{d}) = \left\{ f : f = \sum_{i=1}^{n} a_{i} f_{i} ,
\, \, a_{i} \in\mathbb{R}, \,\, n \in\mathbb{N}, \,\, f_{i} \in C_{c}
(\mathbb{R}^{d}) \cup I(\mathbb{R}^{d}), \right\} \end{align*} We will
sometimes emphasise the dependence on the initial condition by writing
$X_{T} = X_{T}^{x}$, $\mathbb{E}= \mathbb{E}^{x}$, etc.. The following
lemma will also be useful 

\begin{axiom} For every $y \in \mathbb{R}$,
$T>0$ and under the assumption $2\kappa \theta > \eta$ the following is true 
\begin{align*} \lim_{\epsilon \downarrow 0} \sup_{x \in \mathbb{R}}
\mathbb{P}(X_{T}^{x} \in (y -\epsilon, y + \epsilon)) = 0. \end{align*}
\end{axiom} \begin{proof} We write $J_{t} = \int_{0}^{t} \int_{E} y (\mu -
\nu)(dy,ds)$ and observe that the distribution of $X_{T}^{x}$ conditional
on $J_{T}$ and $\{W_{t} : 0\leq t\leq T \}$ is Gaussian. Indeed we have
\begin{align*} X_{T}^{x}|_{J_{T}, \{W_{t} : 0\leq t\leq T \}}
\sim  N ( x + \alpha , \beta^{2}), \end{align*} where
\begin{align*} \alpha & = \int_{0}^{T} \left( r - \frac{1}{2}
\sigma_{t}^{2} \right) dt + \rho \int_{0}^{T} \sigma_{t} dW_{t} + J_{T} \\
\beta^{2} & = (1 - \rho^{2}) \int_{0}^{T} \sigma_{t}^{2} dt. \end{align*}
This gives \begin{align*} \sup_{x \in \mathbb{R}^{d}} \mathbb{P}(X_{T}^{x}
\in (y - \epsilon, y + \epsilon) ) & =\sup_{x \in \mathbb{R}} \mathbb{E}[
\mathbb{E}[1_{\{X_{T}^{x} \in (y - \epsilon, y + \epsilon) \}} | J_{T},
\{W_{t} : 0 \leq t \leq T \}]] \\ & = \sup_{x \in \mathbb{R}}
\mathbb{E}\left[ \int_{y- \epsilon}^{y + \epsilon} \frac{1}{\sqrt{2 \pi}
\beta} \exp \left( - \frac{(z - \alpha - x)^{2}}{2 \beta^{2}} \right) dz
\right] \\ & \leq C \epsilon \mathbb{E}[\beta^{-1}] \end{align*} for some
constant $C < \infty$. The proof will be complete if we can show
$\mathbb{E}[\beta^{-1}] < \infty$, but this is true since the
Cauchy-Schwarz inequality gives \begin{align*} \beta^{-1} \leq (T
\sqrt{1-\rho^{2}})^{-1} \left( \int_{0}^{T} \sigma_{t}^{-2}dt \right)^{1/2}
\, \, \, \, \, \, \, \, \, \, \, a.s. \end{align*} Then, from the previous
lemma, \begin{align*} \mathbb{E} \left[ \left( \int_{0}^{T} \sigma_{t}^{-2}
dt \right)^{1/2} \right] & \leq \mathbb{E}\left[ \int_{0}^{T}
\sigma_{t}^{-2} dt \right]^{1/2} \\ & \leq ( T \sup_{0\leq t\leq T}
\mathbb{E}[\sigma_{t}^{-2}] )^{1/2} < \infty. \end{align*} \end{proof}

An application of the extended Bismut-Elworthy-Li formula will give the
following result. \begin{theorem} Suppose that the parameters of the SVJ
model satisfy $2 \kappa \theta > \eta^{2}$ and $f \in
\mathcal{J}(\mathbb{R})$ then, provided $|\rho| < 1$, the following is true
\begin{align*} \frac{\partial}{\partial S_{0}} \mathbb{E}[f(S_{T})] =
\mathbb{E} \left[ f(S_{T}) \int_{0}^{T} \frac{1}{T S_{0}\sqrt{1-\rho^{2}}
\sigma_{s}} dZ_{s} \right]. \end{align*} \end{theorem} \begin{subsubaxiom}
For the purposes of Monte Carlo applications one would make use of the
localised Malliavin technique described in Fourni\'e et al $[14]$, and it is
clear that the class of functions $\mathcal{J}(\mathbb{R})$ is sufficiently
rich for this purpose. In particular, in enables us to deal with digital
payoffs and European call and put option payoffs. \end{subsubaxiom}

\begin{proof} \emph{Step 1} We assume that $f\in C_{c}^{2}(\mathbb{R})
\subset C_{b}^{2}(\mathbb{R})$ and note that this implies $f \in
C_{b}^{2}(\mathbb{R})$ and we let $D\subset \mathbb{R}$ be some arbitrary
compact subset with $x \in D$. It suffices to derive a representation for
$X_{T}$ for $ f \in \mathcal{J}(\mathbb{R})$, the conclusion for $S_{T}$
will then follow by applying the result for $X_{T}$ to the function $f
\circ \exp \in \mathcal{J}(\mathbb{R})$ , and changing the variable of
differentiation to $S_{0}$. \\
 \emph{Step 2} We construct an approximating
sequence of SDEs with solution $X^{N}$ such that $X^{N}\rightarrow X$ a.s.
and such that the extended Bismut-Elworthy-Li formula can be applied for
each $X^{N}$ . To this end we define for $N \geq 2$ \begin{equation}
\begin{split} \begin{pmatrix} X_{t}^{N} \\ \sigma_{t}^{N} \end{pmatrix} = &
\begin{pmatrix} x \\ \sigma_{0} \end{pmatrix} + \int_{0}^{t}
\begin{pmatrix} r - h^{N}(\sigma_{s}^{N}) \\ g^{N}(\sigma_{t}^{N}) -
\frac{\kappa}{2} \sigma_{s}^{N} \end{pmatrix} ds + \int_{0}^{t}
\begin{pmatrix}\sqrt{1-\rho^{2}} p^{N}(\sigma_{s}^{N}) \\ 0 \end{pmatrix}
dZ_{s} \\ & + \int_{0}^{t} \begin{pmatrix} \rho p^{N}(\sigma_{s}^{N}) \\
\frac{\eta}{2} \end{pmatrix} dW_{s} + \int _{0}^{t} \int_{E}
\begin{pmatrix} y \\ 0 \end{pmatrix} (\mu - \nu) (dy,ds) \end{split}
\end{equation} 
where the functions $h^{N}, g^{N}, p^{N}\in
C^{2}_{b}(\mathbb{R})$ are such that
 \begin{equation*} h^{N}(x) =
\begin{cases} \frac{1}{2} x^{2} & \text{ if $|x| \leq N $} \\ 0 & \text{if
$|x| \geq N+1$} \end{cases}
 \end{equation*} 
\begin{equation*} g^{N}(x) =
\begin{cases} \left(\frac{\kappa \theta}{2} - \frac{\eta^{2}}{8} \right)
\frac{1}{x} & \text{if $x \geq \frac{1}{N} $} \\ 0 & \text{ if $x \leq
\frac{1}{2N} $} \end{cases} 
\end{equation*} and
 \begin{equation*} p^{N}(x)
= \begin{cases} \frac{1}{N^{\xi}} & \text{ if $x \leq 0 $} \\ x & \text{if $x \geq
\frac{1}{N^{\xi }}$} \end{cases} 
\end{equation*} where, $\xi = \frac{1}{2}(\frac{\delta}{2} -1)$ and,
as in Lemma $1$, $\delta = \frac{4\kappa \theta}{\eta^{2}}$.  Moreover,  for each $N$, $h^{N}(x)
\leq \frac{1}{2} x^{2}$ for all $x \in \mathbb{R}$,  $g^{N}(x) \leq \left(\frac{\kappa \theta}{2} -
\frac{\eta^{2}}{8} \right) \frac{1}{x}$ for all $x \in [0,\infty)$ and 
 $\frac{1}{2N^{\xi }} \vee x \leq p^{N}(x) \leq 1$ for all $x \in \left[0, \frac{1}{N} \right]$ (similar approximating sequences for the volatility have
been discussed in Ewald $[12]$). 
Next, we
define the stopping times \begin{align*} \tau_{N} = \inf \left\{ t\geq 0 :
\sigma_{t} \leq \frac{1}{N} \right\}, && \zeta_{N} = \inf \left\{ t\geq 0 :
\sigma_{t} \geq N \right\}. \end{align*} Then, it is well known that for
$\eta^{2} < 2 \kappa \theta$ the volatility never hits zero so we have
$\tau_{N}\rightarrow \infty$ a.s. as $N\rightarrow \infty$, and since the
solution to $(19)$ is non-explosive we also have $\zeta_{N} \rightarrow
\infty$ a.s. as $N\rightarrow \infty$. Consequently, for each $t \in
[0,T]$, $X_{t}^{N} = X_{t}$ a.s. on the set $\{\tau_{N^{\xi}} > t, \zeta_{N} > t
\}$ and so $X_{t}^{N} \rightarrow X_{t}$ a.s. as $N\rightarrow \infty$. \\
\emph{Step 3} We confirm that the extended Bismut-Elworthy-Li formula
applies for each $N$ to deduce \begin{align*} \frac{\partial}{\partial
x}\mathbb{E}[f(X_{T}^{N})] = \mathbb{E} \left[ f(X_{T}^{N}) \int_{0}^{T}
\frac{1}{T \sqrt{1-\rho^{2}} p^{N}(\sigma_{s}^{N})} dZ_{s} \right].
\end{align*} The vector fields driving the SDE defining $X^{N}$ satisfy the
conditions of Theorem $2$ so we need only verify that the process
\begin{align*} K_{t}^{N} : = \int_{0}^{t} R(s,x_{s}) \frac{\partial
x_{t}^{z}}{\partial z_{1}} dW_{s} = \int_{0}^{t} \frac{1}{T
\sqrt{1-\rho^{2}} p^{N}( \sigma_{s}^{N})} dZ_{s} \end{align*} is a
martingale for all $N$. But this is immediate from the fact that the
integrand is bounded (by $2N^{\xi}/\sqrt{1-\rho^{2}}$). \\
\emph{Step 4} Next we
check that \begin{align*} \frac{\partial}{\partial x}
\mathbb{E}^{x}[f(X_{T})] = \lim_{N\rightarrow
\infty}\frac{\partial}{\partial x} \mathbb{E}^{x}[f(X_{T}^{N})]
\end{align*} To do this we define the sequence of functions $\phi ^{N} : D
\rightarrow \mathbb{R}$ by $\phi ^{N} : x \mapsto
\mathbb{E}^{x}[f(X_{T}^{N})]$. We know that each $\phi^{N}$ is
differentiable, and it is clear by bounded convergence that $\phi^{N}(x)
\rightarrow \phi(x) : = \mathbb{E}^{x}[f(X_{T})]$ for every $x \in D$. We
now confirm that $\phi$ is differentiable with $\phi'(x) =
\lim_{N\rightarrow \infty} (\phi^{N})'(x) = \mathbb{E}[f'(X_{T})]$. Since,
for every $N$, $\frac{\partial X_{T}^{N}}{\partial x} \equiv 1$ this will
be achieved if we can show \begin{align} \lim_{N \rightarrow \infty}
\mathbb{E}[f'(X_{T}^{N})] = \mathbb{E}[f'(X_{T})] \end{align} and the
convergence is uniform over $x \in D$. To do this we first show that
$X_{T}^{N} \rightarrow X_{T}$ in $L^{1}$ uniformly in $x \in D$, but since
each term $\mathbb{E}[|X_{T} - X_{T}^{N}|]$ is independent of $x$ it
suffices the show that $X_{T}^{N} \rightarrow X_{T}$ in $L^{1}$, since any
convergence will then immediately be uniform in $x$. Before we do this we
note that a straight forward application of the comparison theorem (page
269 Rogers and Williams $[24]$) tells us for each $t \in [0,T]$ that
$y_{t}\leq \sigma_{t}^{N}$ a.s. where $y_{t}$ is the Ornstein-Uhlenbeck
process solving the SDE \begin{align} dy_{t} = -\frac{\kappa}{2}y_{t} dt +
\frac{\eta}{2} dW_{t}. \end{align} We may also use the proof of the
comparison theorem combined with the fact that $\sigma_{t}>0$ a.s. to show
that $\sigma_{t}^{N}\leq \sigma_{t}$ a.s. for each $t\in[0,T]$. Then, we
use the Cauchy-Schwarz and Burkholder-Davis-Gundy inequalities together
with the fact that $ \left(\sum_{i=1}^{k} x_{i} \right)^{2} \leq k
\sum_{i=1}^{k} x_{i}^{2}$ to show that the family $\{X_{T}^{N} : N \geq 2
\}$ is bounded in $L^{2}$. We end up with \begin{align*}
\mathbb{E}[(X_{T}^{N})^{2}] \leq 4x^{2} + \int_{0}^{T}\left( 8T r^{2} + 8T
\mathbb{E}[h^{N}(\sigma_{t}^{N})^{2}] + 4
\mathbb{E}[p^{N}(\sigma_{t}^{N})^{2}] \right)dt \end{align*} Since
$h^{N}(x) \leq x^{2}/2$, $p^{N}(x) \leq 1 + |x|$, $\sup_{0\leq t\leq T}
\mathbb{E}[ (\sigma_{t}^{N})^{4}]\leq \sup_{0\leq t\leq T}
\mathbb{E}[\sigma_{t}^{4}] + \sup_{0\leq t\leq T} \mathbb{E}[y_{t}^{4}] <
\infty$ we conclude \begin{align*} \mathbb{E}[(X_{T}^{N})^{2}] \leq 4x^{2}
+8r^{2}T^{2} + 2T^{2} \sup_{0\leq t\leq T} \mathbb{E}[(\sigma_{t}^{N})^{4}]
+ 8T( 1+ \sup_{0\leq t\leq T}\mathbb{E}[(\sigma_{t}^{N})^{2}]) \end{align*}
and the right hand side of the inequality may be bounded uniformly in $N$,
and consequently $X_{T}^{N} \rightarrow X_{T}$ in $L^{1}$ uniformly in $x$.
Finally, we verify $(22)$ by noting that \begin{align*}
| \mathbb{E}[f'(X_{T}^{N}) - f'(X_{T})] | \leq \mathbb{E}[ &
| | f'(X_{T}^{N})- f'(X_{T})| 1_{ \{|X_{T}^{N}-X_{T}| \leq \epsilon \}}] \\
& + \mathbb{E}[ |f'(X_{T}^{N})-f'(X_{T})| 1_{ \{|X_{T}^{N}-X_{T}| >
\epsilon \}}] \end{align*} The first term on the right converging to zero
uniformly in $x$ by the uniform continuity of $f'$ and the second term
likewise by the convergence in probability (from Chebyshev's inequality) of
$X_{T}^{N}$ to $X_{T}$ uniformly for $x \in D$ and the boundedness of $f'$.
\\
\emph{Step 5} We now establish \begin{align*} \mathbb{E} \left[ f(X_{T})
\int_{0}^{T} \frac{1}{T \sqrt{1-\rho^{2}} \sigma_{s}} dZ_{s} \right] =
\lim_{N \rightarrow \infty} \mathbb{E} \left[ f(X_{T}^{N}) \int_{0}^{T}
\frac{1}{T \sqrt{1-\rho^{2}} p^{N}(\sigma_{s}^{N})} dZ_{s} \right]
\end{align*} from which it follows that \begin{align}
\frac{\partial}{\partial x} \mathbb{E}[f(X_{T})] = \mathbb{E} \left[
f(X_{T}) \int_{0}^{T} \frac{1}{T \sqrt{1-\rho^{2}} \sigma_{s}} dZ_{s}
\right] . \end{align} To do this we show that $K_{T}^{N} \rightarrow K_{T}$
in $L^{1}$, where 
\begin{align*} K_{T}^{N} = \int_{0}^{T}\frac{1}{ T
\sqrt{1- \rho^{2}}p^{N} (\sigma_{s}^{N})} dZ_{s}, \,\,\, K_{T} = \int_{0}^{T}
\frac{1}{T \sqrt{1- \rho^{2}} \sigma_{s}} dZ_{s}, \end{align*}
from which the required result follows immediately by the boundedness of $f$. 
As a preliminary to this we note that by the time changed squared Bessel representation used in Lemma $1$ we can deduce
\begin{align*}
\mathbb{P}(\tau_{N} \leq T )&  = \mathbb{P} \left( \inf_{0\leq t\leq T} \sigma_{t}^{2} \leq \frac{1}{N} \right) 
 \leq \mathbb{P}\left(\inf_{0\leq t\leq T} Y\left(\frac{\eta^{2}}{4 \kappa} (e^{\kappa t} -1) \right) \leq \frac{e^{\kappa T}}{N} \right) \\
& \leq \mathbb{P} \left( \inf_{0\leq t < \infty} Y_{t} \leq \frac{e^{\kappa T}}{N} \right) = \frac{e^{2\kappa \xi T}}{ N^{2\xi}}.
\end{align*}
The last line following from the observation that the scale function for a $\delta$-dimensional squared Bessel process
is $s(x) = -x^{-2\xi}$ (see page 286 of Rogers and Williams $[24]$). We now observe 
 by  the It\^o-isometry that 
 \begin{align*}
 \mathbb{E}[|K_{T}^{N}&- K_{T}|^{2}]  = \frac{1}{T (1- \rho^{2})} \mathbb{E} \left[
\int_{0}^{T}\left( \frac{1}{\sigma_{s}} - \frac{1}{p^{N}(\sigma_{s}^{N})} \right)^{2}
ds \right] \\
& \leq \frac{2}{T(1 - \rho^{2})} \mathbb{E}\left[
\int_{0}^{T} \left( \left(\frac{1}{\sigma_{s}} - \frac{1}{p^{N}(\sigma_{s})} \right)^{2} + \left( \frac{1}{p^{N}(\sigma_{s})}-
\frac{1}{p^{N}(\sigma^{N}_{s})}\right)^{2} \right) ds \right] .
\end{align*} 
Using the three facts  $p^{N}(x) \geq x$ for all $x$, $p^{N}(\sigma_{t}^{N}) = p^{N}(\sigma_{t})$ for $t< \tau_{N}$ and
$p^{N}(x) \geq  \frac{1}{2N^{\xi}}$ for all $x$, we see that 
\begin{align*}
\mathbb{E}[|K_{T}^{N} -  K_{T}|^{2}]   \leq  & \frac{4 \sup_{0\leq t\leq T} \mathbb{E}[ \sigma_{t}^{-2}]}{(1 - \rho^{2})} \\
& + \frac{2}{T(1-\rho^{2})} \mathbb{E}\left[ \int_{\tau_{N}}^{T}  \left( \frac{1}{p^{N}
(\sigma_{s})} - \frac{1}{p^{N}(\sigma_{s}^{N})} \right)^{2}dt \, 1_{\{\tau_{N} < T \}} \right] \\
& \leq \frac{4}{1-\rho^{2}} \left( \sup_{0\leq t\leq T} \mathbb{E}[\sigma_{t}^{-2}] + 4N^{2\xi} \mathbb{P}(\tau_{N} < T) \right) \\
& \leq \frac{4}{1-\rho^{2}} \left( \sup_{0\leq t\leq T} \mathbb{E}[\sigma_{t}^{-2}] + 4e^{2 \xi \kappa T} \right) < \infty,
\end{align*}
and the fact that $K_{T}^{N} \rightarrow K_{T}$ in $L^{1}$ is immediate. \\
\emph{Step 6} We now relax the
regularity conditions on $f$ in two stages. Firstly, we extend to $f \in
C_{c}(\mathbb{R})$. To do this we notice that we can identify a sequence of
functions $f_{n} \in C_{c}^{\infty}(\mathbb{R})$ with $f_{n}\rightarrow f$
uniformly and boundedly as $n\rightarrow \infty$. The extension is then
immediate since bounded convergence implies $\mathbb{E}[f_{n}(X_{T})]
\rightarrow \mathbb{E}[f(X_{T})]$ and, for any compact subset $H \subset
\mathbb{R}$, we have \begin{align} \sup_{x \in H} \left|
\frac{\partial}{\partial x} \mathbb{E}[f_{n}(X_{T}^{x})] -
\mathbb{E}[f(X_{T}) K_{T}] \right| \leq \mathbb{E}[K_{T}^{2}]^{1/2}\sup_{x
\in H} \mathbb{E}[(f_{n}(X_{T}^{x}) - f(X_{T}^{x}))^{2}]^{1/2} \end{align}
The convergence of the right hand side to zero being immediate from the
fact that $f_{n} \rightarrow f$ uniformly, and that $X_{T}^{x} = x + S$ for
some random variable $S$ independent of $x$. Secondly, we extend to
indicator functions of the form $f = 1_{[a,b]}$ (the extension to
indicators of open and half-open intervals being similar). To do this we
note that we can construct an approximating sequence $f_{n} \in
C_{c}(\mathbb{R})$ having the properties that $f_{n} \rightarrow f$
pointwise and, for any neighbourhoods $B_{a}$ and $B_{b}$ of $a$ and $b$
respectively, $f_{n}-f = 0$ on $L : = \mathbb{R}\cap B_{a}^{c} \cap
B_{b}^{c}$ for $n$ sufficiently large. We can now repeat the argument of
the previous paragraph to obtain $(25)$. To show that the right hand side
of $(25)$ can be made arbitrarily small we let $\delta > 0$ and fix some
$\epsilon > 0$ chosen such that \begin{align*} \sup_{x \in H}
\mathbb{P}(X_{T}^{x} \in (a -\epsilon, a + \epsilon)) < \frac{\delta}{4
\mathbb{E}[K_{T}^{2}]^{1/2}} \end{align*} and \begin{align*} \sup_{x \in H}
\mathbb{P}(X_{T}^{x} \in (b -\epsilon, b + \epsilon)) < \frac{\delta}{4
\mathbb{E}[K_{T}^{2}]^{1/2}} \end{align*} as we may by Lemma $2$. With $L =
\mathbb{R} \cap (a-\epsilon, a+ \epsilon)^{c} \cap (b-\epsilon, b+
\epsilon)^{c}$ we may then choose $N$ such that for all $n\geq N$ we have
$\sup_{y \in L}|f_{n}(y) - f(y)| = 0 $ and we can bound the right hand side
of $(26)$ by \begin{align*} 2 \mathbb{E}[K_{T}^{2}]^{1/2} \sup_{x \in H}
\big( \mathbb{P}(X_{T}^{x} \in (a -\epsilon, a+ \epsilon)) +
\mathbb{P}(X_{T}^{x} \in (b -\epsilon, b+ \epsilon)) \big) < \delta.
\end{align*} Since $\delta$ was arbitrary this completes the result. Since
it is clear that $(24)$ is stable under taking finite linear combinations
the extension to the class $\mathcal{J}(\mathbb{R})$ is immediate. The
result for $S_{T}$ follows as described in Step $1$. \end{proof}
\begin{subsubaxiom} By the same argument and under the same conditions as
the last theroem we can also obtain \begin{align*} \frac{\partial}{\partial
S_{0}} \mathbb{E}[f(S_{T_{1}},\ldots,S_{T_{n}})] = \mathbb{E} \left[
f(S_{T_{1}},\ldots,S_{T_{n}}) \int_{0}^{T_{1}} \frac{1}{T
S_{0}\sqrt{1-\rho^{2}} \sigma_{s}} dZ_{s} \right], \end{align*} for any $n
\in \mathbb{N}$ and $0<T_{1}\leq T_{2}\leq \ldots \leq T_{n}\leq T$.
\end{subsubaxiom}

\begin{subsubaxiom} We may apply Theorem $3$ together with a similar
approximation procedure described above to deduce the representation for
the gamma \begin{align*} & \frac{\partial^{2}}{\partial S_{0}^{2}}
\mathbb{E}[f(S_{T})] = \\ &\frac{4}{(1-\rho^{2})T^{2}S_{0}^{2}}
\mathbb{E}\left[\left( \int_{T/2}^{T} \frac{1}{\sigma_{t}}dZ_{t}
\int_{0}^{T/2} \frac{1}{\sigma_{t}}dZ_{t} - \frac{T \sqrt{1-\rho^{2}}}{4}
\int_{0}^{T} \frac{1}{\sigma_{t}}dZ_{s} \right) f(S_{T}) \right] \\ & =
\mathbb{E} \left[ \frac{4}{(1-\rho^{2})T^{2}S_{0}^{2}}\left( \int_{T/2}^{T}
\frac{1}{\sigma_{t}}dZ_{t} \int_{0}^{T/2} \frac{1}{\sigma_{t}}dZ_{t}\right)
f(S_{T}) \right] - \frac{1}{S_{0}} \frac{\partial}{\partial S_{0}}
\mathbb{E}[f(S_{T})] \end{align*} for $f \in \mathcal{J}(\mathbb{R})$.
Where, as above, we have initially used Theorem $3$ for $X_{T}$ and deduced
the result for $S_{T}$ by applying it to the function $f \circ \exp $ and
using the observation that \begin{align*} \frac{\partial^{2}}{\partial
x^{2}} = S_{0} \frac{\partial}{\partial S_{0}} + S_{0}^{2}
\frac{\partial^{2}}{\partial S_{0}^{2}} \end{align*} \end{subsubaxiom}

\subsection{Stochastic volatility with jumps in the volatility - the
Matytsin model}

We consider how these ideas may be extended to the model of Matytsin where
the volatility evolves according to the Heston model with the exception
that there are jumps which occur in the stock and volatility
simultaneously, the volatility jumps being of positive deterministic size.
This volatility process is written as \begin{align*} d\sigma_{t}^{2} =
\kappa( \theta- \sigma_{t}^{2}) dt + \eta\sigma_{t}dW_{t} + \gamma dJ_{t}
\end{align*} where $J_{t}$ is a Poisson process. With $X_{t}$ given as in
the the SVJ model and $\sigma_{t}^{2}$ as above the pair
$(X_{t},\sigma_{t}^{2})$ describes the Matytsin double jump model ( or
SVJJ). Applying It\^o's formula we can express the system
$(X_{t},\sigma_{t})$ in terms of our previous notation by the SDEs \[
\begin{split} \begin{pmatrix} X_{t}\\ \sigma_{t} \end{pmatrix} =
\begin{pmatrix} x\\ \sigma_{0} \end{pmatrix} & + \int_{0}^{t}
\begin{pmatrix} r - \frac{1}{2}\sigma_{s}^{2}\\ \left(
\frac{\kappa\theta}{2} - \frac{\eta^{2}}{8} \right) \frac{1} {\sigma_{s}}
- \frac{\kappa}{2} \sigma_{s} + \lambda( \sqrt{\sigma_{s-}^{2} + \gamma} -
\sigma_{s-}) \end{pmatrix} ds\\ & + \int_{0}^{t} \begin{pmatrix}
\rho\sigma_{s}\\ 0 \end{pmatrix} dZ_{s} + \int_{0}^{t} \begin{pmatrix}
\sqrt{1-\rho^{2}}\sigma_{s}\\ \frac{\eta}{2} \end{pmatrix} dW_{s}\\ & +
\int_{0}^{t} \int_{E} \begin{pmatrix} y\\ \sqrt{\sigma_{s-}^{2} +\gamma} -
\sigma_{s-} \end{pmatrix} (\mu- \nu) (dy,ds) \end{split} \] where $E =
\mathbb{R}$, $\gamma$ is the constant jump size in the volatility and $\mu$
is a Poisson random measure with mean measure $v(dy,dt) = \lambda G(dy)dt =
\lambda p(y) dy dt$ where here $p(y)$ is the density of the jumps in $X$.

\begin{theorem} Suppose that the parameters in the SVJJ model satisfy
$2\kappa \theta > \eta^{2}$ and $f \in \mathcal{J}(\mathbb{R})$. Then,
provided $|\rho|< 1$, the following is true \begin{align*}
\frac{\partial}{\partial S_{0}} \mathbb{E}[f(S_{T})] = \mathbb{E} \left[
f(S_{T}) \int_{0}^{T} \frac{1}{ S_{0}T \sqrt{1-\rho^{2}} \sigma_{s-}}
dZ_{s} \right] \end{align*} \end{theorem} \begin{proof} The proof may be
completed by following the steps of the previous theorem. The approximating
system used is \begin{equation*} \begin{split} \begin{pmatrix} X_{t}^{N} \\
\sigma_{t}^{N} \end{pmatrix} = \begin{pmatrix} x \\ \sigma_{0}
\end{pmatrix} & + \int_{0}^{t} \begin{pmatrix} r- h^{N}(\sigma_{s-}^{N}) \\
g^{N}(\sigma_{s-}^{N}) - \frac{\kappa}{2} \sigma_{s}^{N} + \lambda (
\sqrt{(\sigma_{s-}^{N})^{2} + \gamma} - \sigma_{s-}^{N}) \end{pmatrix} ds
\\ & + \int_{0}^{t} \begin{pmatrix}\rho p^{N}(\sigma_{s-}^{N}) \\ 0
\end{pmatrix} dZ_{s} + \int_{0}^{t} \begin{pmatrix}
\sqrt{1-\rho^{2}}p^{N}(\sigma_{s-}^{N}) \\ \frac{\eta}{2} \end{pmatrix}
dW_{s} \\ & + \int _{0}^{t} \int_{E} \begin{pmatrix} y \\
\sqrt{(\sigma_{s-}^{N})^{2} + \gamma} - \sigma_{s-}^{N} \end{pmatrix} (\mu
- \nu) (dy,ds) \end{split} \end{equation*} and then, by the same argument
as before, $\sigma_{t}^{N} \rightarrow \sigma_{t}$ almost surely. Denoting
$\tilde{\sigma}$ to be the solution of the usual continuous Heston process
with the same parameters and its approximating process by
$\tilde{\sigma}^{N}$, and using the fact that the jumps in the volatility
in Matytsin are non-negative we can apply the comparison theorem in between
jumps to give the relation \begin{align} y_{t} \leq \tilde{\sigma}^{N}_{t}
\leq \sigma_{t}^{N} \leq \sigma_{t} \end{align} a.s. for every $t \in
[0,T]$, where $y_{t}$ is as in $(23)$. Consequently, using $(x + y)^{p}
\leq 2^{p-1}(x^{p} + y^{p})$ we can deduce \begin{align*}
\mathbb{E}[\sup_{0\leq t \leq T} (\sigma_{s}^{N})^{4} ] \leq 8
\mathbb{E}[\sup_{0\leq t\leq T} \sigma_{t}^{4}] + 8 \mathbb{E}[\sup_{0\leq
t\leq T} y_{t}^{4}] < \infty \end{align*} and so $\sigma_{t}^{N}
\rightarrow \sigma_{t}$ in $L^{2}$ for every $t \in [0,T]$ by dominated
convergence. The remainder of the argument follows as before, only the
elementary observation (which follows from $ (26)$) that $\sup_{0\leq t\leq
T} \mathbb{E}[\sigma_{t}^{-2}] \leq \sup_{0\leq t\leq T}
\mathbb{E}[\tilde{\sigma}_{t}^{-2}]$ is needed to recycle the estimates
already established for the process $\tilde{\sigma }$ to give new estimates
on the Matytsin volatility $\sigma$. The extension from
$C_{b}^{2}(\mathbb{R})$ to $\mathcal{J}(\mathbb{R})$ proceeds in the same
way as before after Lemma $2$ has been verified with the Matytsin
volatility, which follows from an elementary adaptation of the argument
given. \end{proof}

\begin{subsubaxiom} Under the same assumptions of Remark $8$ we may again
derive an representation for the gamma for the SVJJ model analogous to the
one for SVJ. \end{subsubaxiom}

\section{Numerical Results}

We implement the results for the SVJJ model firstly in the case of a
European call option (payoff $(S_{T}-K)_{+}$) with $T=1$ and $S_{0} = 100$
and strike $K = 100$, secondly for a double digital payoff of the form
$1_{[K_{1},K_{2} ]}(S_{T})$ again with $T=1$,$S_{0} =100$ and $K_{1} =
100$, $K_{2} =110$. Finally, we implement for the delta of a digital
Cliquet option with payoff profile
$1_{[K^{*}_{1},K^{*}_{2}]}(S_{T}-S_{T_{1}})$, where $T=1,T_{1}=0.5,
K_{1}^{*}=5$ and $K_{2}^{*}=10$. The model parameters we use are $r=0$,
$\rho= -0.7$, $\gamma=0.4$, $\sigma_{0}^{2} = 0.1$, $\lambda= 1.0$,
$\theta= 0.08$, $\kappa=4.0$, $\eta=0.6$, and we assume that the jumps in
$\log S$ are distributed normally with mean $-0.1$ and standard deviation
$0.1$. 
\newline
\newline
\newline
\newline
\newline
\newline
\begin{figure}[H]
\scalebox{1.4}{\includegraphics[height = 6cm, width
=8cm]{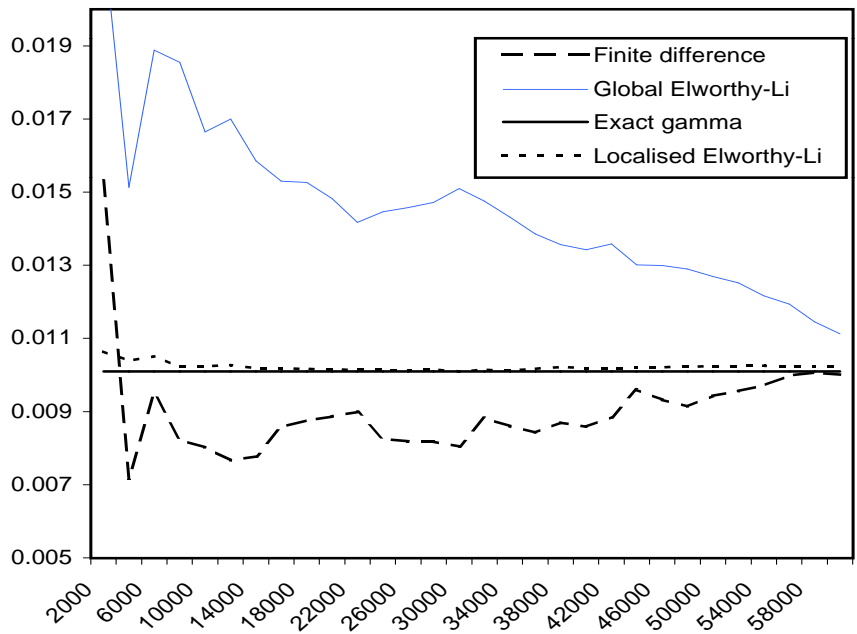}}\caption{Gamma for a European call option with parameters
as above} \end{figure}
\begin{figure}[H]
\scalebox{0.9}{\includegraphics[height =8cm, width
=12cm]{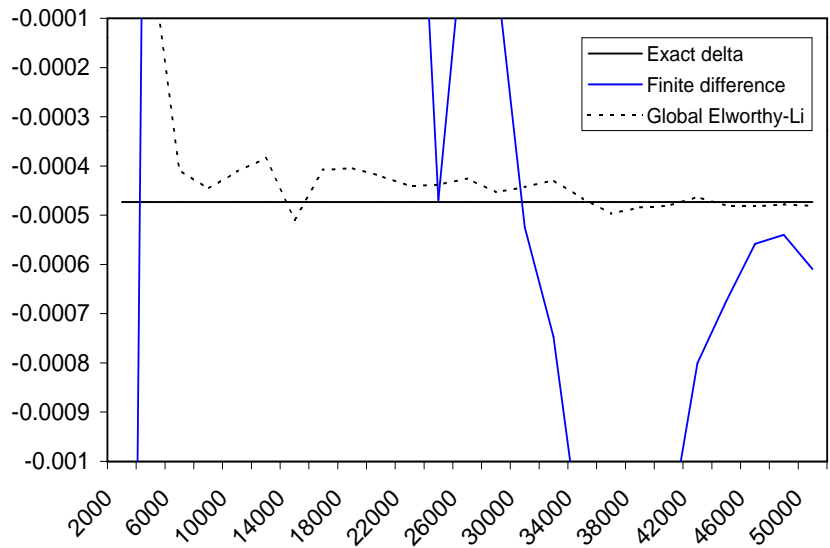}}\caption{Delta for a double digital option with
parameters as above} \end{figure}
\begin{figure}[H]
\scalebox{0.9}{\includegraphics[height =8cm, width
=12cm]{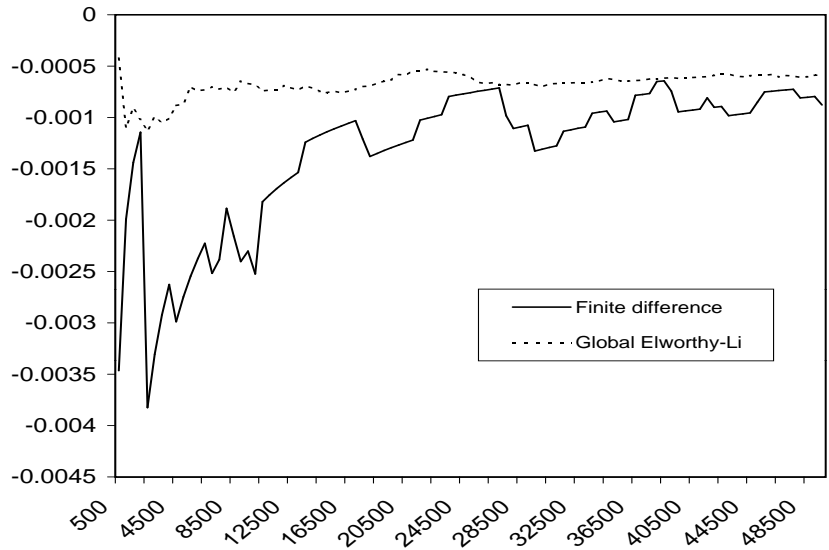}}\caption{Delta for a digital Cliquet option with parameters
as above} 
\end{figure}

\end{document}